\documentclass{article}[11]

\usepackage{latexsym}
\usepackage{amssymb} 
\usepackage{epsfig}
\usepackage{amsmath}
\usepackage[mathscr]{eucal}

\newtheorem{Lemma}{Lemma}
\newtheorem{Proposition}[Lemma]{Proposition}
\newtheorem{Theorem}[Lemma]{Theorem}
\newtheorem{Example}{Example}

\newtheorem{Definition}{Definition}

\newtheorem{Conjecture}{Conjecture}

\newtheorem{Question}{Question}

\newcommand{\ed}{\ \stackrel{d}{=} \ }

\newcommand{\EE}{\mbox{${\mathcal E}$}}
\newcommand{\VV}{\mbox{${\mathcal V}$}}

\newcommand{\SSS}{\mbox{${\mathcal S}$}}

\newcommand{\eps}{\varepsilon}
\newcommand{\bE}{{\bf E}}
\newcommand{\bP}{{\bf P}}

\newcommand{\bx}{{\bf x}}

\newcommand{\bzero}{{\bf 0}}

\newcommand{\MMM}{\mbox{${\mathfrak{M}}$}}

\newcommand{\Pbb}{\mbox{${\mathbb P}$}}

\newcommand{\sfrac}[2]{{\textstyle\frac{#1}{#2}}}

\newcommand{\qed}{\hfill{\ \ \rule{2mm}{2mm}} \vspace{0.2in}}
\newcommand{\proof}{\noindent \emph{Proof:}\ }

\newcommand{\Tbold}{{\mathbb{T}}}
\newcommand{\Zbold}{{\mathbb{Z}}}
\newcommand{\Nbold}{{\mathbb{N}}}

\newcommand{\logg}{\log_2}
\newcommand{\Aut}{{\rm Aut}}
\newcommand{\dist}{{\rm dist}}

\title{{\bf On the Cluster Size Distribution for 
Percolation on Some General Graphs}} 

\author{{\bf Antar Bandyopadhyay}
\footnote{Theoretical Statistics and Mathematics Unit, 
Indian Statistical Institute, Delhi Centre, 7 S. J. S. Sansanwal Marg, 
New Delhi 110016, India, E-Mail: {\tt antar@isid.ac.in}}
\and 
{\bf Jeffrey Steif}
\footnote{Mathematical Sciences, Chalmers University of Technology
and Mathematical Sciences, G\"{o}teborg University, 
SE-412 96, G\"{o}teborg, Sweden, E-Mail: {\tt steif@math.chalmers.se}}
\and
{\bf \'{A}d\'{a}m Tim\'{a}r}
\footnote{Mathematics Department, The University of British Columbia,
1984 Mathematics Road, Vancouver V6T 1Z2, Canada, 
E-Mail: {\tt timar@math.ubc.ca}}}

\begin{document}

\maketitle

\begin{abstract}
We show that for any Cayley graph, the probability (at any $p$)
that the cluster of the origin has size $n$ decays at a well-defined
exponential rate (possibly $0$). For general graphs, we 
relate this rate being positive in the supercritical regime with the 
amenability/nonamenability of the underlying graph.

\vspace{0.1in}

\vspace{0.1in}
\end{abstract}

\section{Introduction}
\label{Sec:Introduction}
Percolation is perhaps the most widely studied statistical physics model
for modeling random media. In addition,
it is a source of many challenging open problems and beautiful conjectures
which are easy to state but often are very difficult to settle;
see \cite{Gri99} for a survey and introduction. The classical literature
concentrates on studying the model on Euclidean lattices $\Zbold^d$, 
$d \geq 2$ and on trees. However in recent years,
there has been a great deal of interest
in studying percolation on other infinite, locally finite, connected 
graphs; see
\cite{BenSch96, BLSP99a, BLSP99b, HagPeres99, HagPeresSch99, Hag00, 
HagSchStei00, Sch01}.  

Our first theorem states that for any Cayley graph, the probability 
that the cluster of the origin has size $n$ decays at a well-defined
exponential rate. For $\Zbold^d$, this is Theorem 6.75 in \cite{Gri99}.

Throughout this paper, $C$ will denote the connected component of
a fixed vertex (the origin for Cayley graphs) for Bernoulli percolation.

\begin{Theorem} 
\label{Thm:decay}
If $G$ is a Cayley graph, then
\[ \lim_{n \rightarrow \infty} -{1\over n} \log\, 
\Pbb_p\left(|C|=n \right) \]
exists for every $p \in \left(0,1\right)$. 
\end{Theorem}

Our method for proving this result combines a 
randomized version of the usual method using subadditivity 
(as in for $\Zbold^d$)
together with a proof that any two finite subgraphs of
$G$ have disjoint translates that are at distance $\leq \delta$ from each
other where $\delta$ is an appropriate function of the sizes of the
subgraphs. One expects perhaps that one should be able to take
$\delta$ being a constant, depending only on the graph. See
Question \ref{Q:touch} for the statement of this problem.

\medskip\noindent
{\bf Remark:} Interesting, as we point out later, there is a concept
of an ordered group and for such groups, the proof of
Theorem 6.75 in \cite{Gri99} can be extended. However, for general groups,
it seems that this proof cannot be applied.

\medskip\noindent
It is of course of interest to know if the limit above is positive or 0.
As will be pointed out later, it is positive below the critical value
for all transitive graphs and so we restrict discussion to
the supercritical regime. In this case, for $\Zbold^d$, the limit is 0
(see Theorem 8.61 in \cite{Gri99}) while for trees it is positive
(although 0 at the critical value). 
Equation (10.12) in \cite{Gri99}) has an explicit formula
for these probabilities for the rooted infinite $3$-ary tree.

One of the key issues studied in percolation
is the difference in the behavior
of percolation depending on whether the underlying graph is \emph{amenable}
or \emph{nonamenable} \cite{BLSP99a, BLSP99b, HagSchStei00, Sch01}. 
For example, for amenable transitive graphs, there is uniqueness of the
infinite cluster for all values of $p$ while for
nonamenable transitive graphs, it is conjectured that there is 
nonuniqueness of the infinite cluster for some values of $p$.
Here it is also worthwhile to point out that
it is well known that properties of other probabilistic models associated
with a graph differ depending on whether the graph is amenable or not. 
Perhaps the most classical of all is the relation with simple random walk 
on a graph, first studied by Kesten \cite{Kes59} where it was shown that
there is a positive spectral gap in the transition operator if and only if
the group is nonamenable. Similar relationships have been
investigated with respect to other statistical physics models 
(see e.g. \cite{Jon99, JonSti99, HagSchStei00, ChenPeres04}).

For the nonamenable case, we state the following question.

\begin{Question}
\label{Q:Transitive-Nonamenable-Exp-Decay}  
Is it true that for a general transitive nonamenable graph $G$ we have
\begin{equation} 
\Pbb_p\left( \left\vert C \right\vert = n \right) \leq 
\exp\left( - \gamma\left(p\right) n \right) \,\,\,\, \forall\,\,\, n \geq 1 \,
\label{Equ:Transitive-Nonamenable-Exp-Decay}
\end{equation}
for some $\gamma\left(p\right) > 0$ whenever $p \neq p_c\left(G\right)$ ? 
\end{Question}

Consider a
general weakly nonamenable graph $G := \left( \VV, \EE \right)$
(not necessarily transitive) with bounded degree. Using a not so difficult
argument of counting \emph{lattice animals}, one can prove that if
$v_0$ is a fixed vertex of $G$ and $C$ is the open connected
component of $v_0$, then for sufficiently large $p$ there is a function
$\gamma\left(p\right) > 0$, such that
\begin{equation}
\Pbb_p\left( \vert C \vert = n \right) \leq e^{- \gamma\left(p\right) \, n }
\,\,\,\, \forall \,\,\, n \geq 1 \,.
\label{Equ:Exp-decay-nonamenable}
\end{equation}

In fact, in the appendix by G\'{a}bor Pete 
in \cite{ChenPeres04} (see equation (A.3)),
it is shown by a slightly more involved argument, 
that the exponential decay
(\ref{Equ:Exp-decay-nonamenable}) holds whenever
$p > 1/\left(1 + \kappa'\right)$ where
$\kappa' = \kappa'\left( G, v_0 \right)$ is the 
anchored Cheeger constant. 
This is certainly in contrast to
the  $\Zbold^d$ case and also,
as we will see later in Section \ref{Sec:Amenable-SubExp-Decay},
to what happens for a large class of transitive amenable graphs.

Using classical branching process arguments,
one can conclude that for any infinite regular tree (which are prototypes
for transitive nonamenable graphs), we must have an
exponential tail bound for the cluster size distribution, when 
$p$ is not equal to the critical probability. 

The assumption of transitivity is however needed for
Question \ref{Q:Transitive-Nonamenable-Exp-Decay} to have a positive
answer as the following example illustrates.
The graph obtained by taking $\Zbold^d$ and attaching a regular rooted tree
with degree $r+1$ at each vertex with $p_c\left(\Zbold^d\right) < \sfrac{1}{r}$
is a nontransitive, nonamenable graph which possesses an intermediate regime
(above the critical value) of sub-exponential decay as next stated
in detail.

\begin{Theorem}
\label{Thm:Non-Transitive-Two-Regime} 
Consider the graph just described and
suppose $p_c\left(\Zbold^d\right) < \sfrac{1}{r}$.
\begin{itemize}
\item[(a)] If
           $p \in \left(0,p_c\left(\Zbold^d\right)\right) \cup
           \left(\sfrac{1}{r}, 1\right)$ then there are functions
           $\phi_1\left(p\right) < \infty$ and $\phi_2\left(p\right) > 0$,
           such that for all $n \geq 1$,
           \begin{equation}
           \exp\left(- \phi_1\left(p\right) n \right) \leq 
           \Pbb_p\left( n \leq \vert C \vert < \infty \right) 
           \leq \exp\left( - \phi_2\left(p\right) \, n \right).
           \label{Equ:Non-Transitive-Exp-Decay}
           \end{equation}

\item[(b)] If $p \in \left(p_c\left(\Zbold^d\right), \sfrac{1}{r}\right)$
           then there are functions $\psi_1\left(p\right) < \infty$ 
           and $\psi_2\left(p\right) > 0$,
           such that for all $n \geq 1$,
           \begin{equation}
           \exp\left( - \psi_1(p) \, n^{(d-1)/d} \right) \leq 
           \Pbb_p\left( n \leq \vert C \vert < \infty \right) \leq 
           \exp\left( - \psi_2(p) \, n^{(d-1)/d} \right).
           \label{Equ:Non-Transitive-Subexp-Decay}
           \end{equation}

\item[(c)] For $p = \sfrac{1}{r}$ we have constants
 $c_1 > 0$ and $c_2\left(\eps\right) < \infty$ such that
           every $\eps > 0$ and for all $n \geq 1$,
           \begin{equation}
           \frac{c_1}{n^{1/2}} \leq 
           \Pbb_p\left( n \leq \vert C \vert < \infty \right) 
           \leq \frac{c_2}{n^{1/2 - \eps}}. 
           \label{Equ:Non-Transitive-Poly-Decay}
           \end{equation}
    \end{itemize}
\end{Theorem} 

\medskip\noindent
As also explained in Section \ref{Sec:Special-Non-Transitive-Graph},
if $p_c\left(\Zbold^d\right) > \sfrac{1}{r}$, this intermediate 
regime disappears.

An interesting class of graphs to investigate in regard to
Question \ref{Q:Transitive-Nonamenable-Exp-Decay}
are products of $\Zbold^d$ with a homogeneous tree.

\begin{Question}
Is there exponential decay in the supercritical regime for 
$\Zbold^d\times T_r$ where $T_r$ is the homogeneous $r$-ary tree?
\end{Question}

We now move to the amenable case.

\begin{Conjecture}
\label{Conj:Transitive-SubExp-Decay}
Let $G := \left( \VV, \EE \right)$ be a transitive amenable graph.
Then there is a sequence $\alpha_n = o\left(n\right)$, such that for
$p > p_c\left(G\right)$
\begin{equation} 
\Pbb_p\left( n \leq \left\vert C \right\vert < \infty \right) \geq
\exp\left( - \eta\left(p\right)\alpha_n \right) 
\,\,\,\, \forall \,\,\, n \geq 1 \,,
\label{Equ:Transitive-SubExp-Decay}
\end{equation}
where $\eta\left(p\right) < \infty$. 
\end{Conjecture}

It turns out that the argument
of Aizenman, Delyon and Aouillard \cite{AizDeSou80, Gri99} 
for proving this sub-exponential behavior
for $\Zbold^d$ can be successfully carried
out for a large class of transitive amenable graphs.
For  $\Zbold^d$, the sequence $\{\alpha_n\}$ can be taken to be 
$\{n^{\frac{d-1}{d}}\}$.

\begin{Theorem}
\label{Thm:Amenable-SubExp-Decay}
If $G := \left( \VV, \EE \right)$ is a Cayley graph of a 
finitely \emph{presented} 
amenable group with \emph{one end}, then 
there is a sequence $\alpha_n = o\left(n\right)$ such that for
$p > p_c\left(G\right)$, there is $\eta\left(p\right) < \infty$ such that
\begin{equation} 
\Pbb_p\left( n \leq \left\vert C \right\vert < \infty \right) \geq
\exp\left( - \eta\left(p\right)\alpha_n \right) 
\,\,\,\, \forall \,\,\, n \geq 1 \,.
\label{Equ:Transitive-SubExp-Tail-Decay}
\end{equation}
\end{Theorem}

We finally point out that transitivity is a necessary condition in
Conjecture \ref{Conj:Transitive-SubExp-Decay}.

\begin{Proposition}
\label{prop:transneeded}
There is an amenable nontransitive graph for which one has exponential
decay of the cluster size distribution at all $p\neq p_c$.
\end{Proposition}

This paper concerns itself mostly with the supercritical case. It therefore 
seems appropriate to end this introduction 
with a few comments concerning the subcritical case. 
It was shown independently in \cite{M86} and \cite{AB87} 
that for $\Zbold^d$ in the subcritical regime,
the size of the cluster of the origin has a finite expected value.
While it seems that the argument in \cite{M86} does not work for all 
transitive graphs as it seems that it is needed
that the balls in the graph grow slower than $e^{n^\gamma}$ for some $\gamma <1$,
it is stated in \cite{Sch01} that the argument in
\cite{AB87} goes through for any transitive graph. 
Theorem 6.75 in \cite{Gri99} (due to \cite{AN84}) states that for 
$\Zbold^d$, if the expected size of the cluster is finite, then 
exponential decay of the tail of the cluster size follows. As stated in
\cite{AN84}, this result holds quite generally in transitive situations and so,
in combination with the statement in \cite{Sch01} referred to above,
for all transitive graphs, one has exponential decay of the
cluster size in the subcritical regime.

We point out however, not surprisingly, 
that transitivity is again needed here. 
An example of a graph which does not have
exponential decay in (a portion of) 
the subcritical regime is obtained by taking the positive
integers, planting a binary tree of depth $a_k$ (sufficiently large) at $k$
for $k\ge 1$ and also attaching to the origin a graph whose critical value 
is say $3/4$. This graph has $p_c=3/4$ but for some $p<3/4$, exponential
decay fails.

We mention that Questions \ref{Q:touch} and \ref{Q:Foelner}
which appear later on and arise naturally in our study 
could also be of interest to people in geometric group theory.

The rest of the paper is organized as follows.
In Section \ref{SubSec:Definitions}, we provide all the necessary 
definitions and notations. 
In Section \ref{Sec:Ordered-Group}, we prove Theorem \ref{Thm:decay}.
In Section \ref{Sec:Amenable-SubExp-Decay},
we prove Theorem \ref{Thm:Amenable-SubExp-Decay}
and Proposition \ref{prop:transneeded}.
Finally, in Section \ref{Sec:Special-Non-Transitive-Graph},
we prove Theorem \ref{Thm:Non-Transitive-Two-Regime} as well as study
the variant of the example in
Theorem \ref{Thm:Non-Transitive-Two-Regime} obtained by taking
$p_c\left(\Zbold^d\right) > \sfrac{1}{r}$ instead.

\section{Definitions and notations}
\label{SubSec:Definitions}
Let $G = \left(\VV, \EE\right)$ be an infinite, connected graph. 
We will say $G$ is \emph{locally finite} if every vertex has finite
degree. 

The i.i.d.\ Bernoulli
bond percolation with probability $p \in \left[0,1\right]$
on $G$ is a probability measure on $\left\{0, 1\right\}^{\EE}$, such that 
the coordinate variables are i.i.d.\ with $\mbox{Bernoulli}\left(p\right)$
distribution. This measure will be denoted by $\Pbb_p$.  
For a given configuration in $\left\{0, 1\right\}^{\EE}$,
it is customary to say that an edge $e \in \EE$ is 
\emph{open} if it is in state $1$, otherwise it is said to be \emph{closed}. 
Given a configuration, write
$\EE = \EE_o \cup \EE_c$, where
$\EE_o$ is the set of all open edges and $\EE_c$ is the set of all closed
edges. The connected components of the subgraph
$\left(\VV, \EE_o\right)$ are called the 
\emph{open connected components} or \emph{clusters}. 

One of the fundamental quantities in
percolation theory is the critical probability $p_c\left(G\right)$ defined by
\begin{equation}
p_c\left(G\right) := \inf \left\{ p \in \left[0,1\right] \,\Big\vert\,
\Pbb_p\left( \exists \mbox{\ an infinite cluster\ } \right) = 1 \,\right\}
\,.
\label{Equ:Defi-Crticical-Prob}
\end{equation}
The percolation model is said to be \emph{subcritical}, \emph{critical} or
\emph{supercritical} regime depending on whether
 $p < p_c\left(G\right)$, 
$p = p_c\left(G\right)$ or $p > p_c\left(G\right)$ respectively. 

For a fixed vertex $v \in \VV$, let $C\left(v\right)$ be the open connected
component containing the vertex $v$. Let 
\begin{equation}
\theta_G^v\left(p\right) := 
\Pbb_p\left( C\left(v\right) \mbox{\ is infinite\ } \right) \, .
\label{Equ:Defi-Perc-Func-v}
\end{equation}
For a connected graph $G$, it is easy to show that irrespective of the
choice of the vertex $v$
\begin{equation}
p_c\left(G\right) = \inf \left\{ p \in \left[0,1\right]
\,\Big\vert\, \theta_G^v\left(p\right) > 0 
\, \right\} \,.
\end{equation}

\begin{Definition}
\label{Defi:Transitive}
We will say a graph $G = \left(\VV, \EE\right)$ 
is \emph{transitive} if for every pair of vertices
$u$ and $v$ there is an automorphism of $G$, which sends $u$ to $v$. In
other words, a graph $G$ is transitive if its automorphism group
${\tt Aut}\left(G\right)$ acts transitively on $\VV$. 
\end{Definition}
Observe that if $G$ is transitive then we can drop the dependency on the
vertex $v$ in (\ref{Equ:Defi-Perc-Func-v}), and then we can write
$\theta_G\left(p\right) = 
\Pbb_p\left( C\left(v_0\right) \mbox{\ is infinite\ } \right)$
for a fixed vertex $v_0$ of $G$. $\theta_G\left(\cdot\right)$ is 
called the \emph{percolation function} for a transitive graph $G$.

We now give definitions of some of the qualitative properties of a
graph $G$ which are important for our study. 

\begin{Definition}
\label{Defi:Amenable}
Let $G := \left(\VV, \EE\right)$ be an infinite, locally finite, connected
graph. The \emph{Cheeger constant} of $G$, denoted by $\kappa\left(G\right)$,
is defined by 
\begin{equation}
\kappa\left(G\right) := \inf \, 
\left\{ \left. \frac{\left\vert \partial W \right\vert}
             {\left\vert W \right\vert}
\,\right\vert\, 
\emptyset \neq W \subseteq \VV \mbox{\ and\ }
\left\vert W  \right\vert  < \infty \, \right\}
\label{Equ:Cheeger-Const}
\end{equation}
where 
$\partial W := \left\{ u \not\in W \,\Big\vert\, \exists v \in 
W, \mbox{\ such that\ } \left\{u, v\right\} \in \EE \,\right\}$
is the \emph{external vertex boundary}.  
The graph $G$ is said to be \emph{amenable} if
$\kappa\left(G\right) = 0$; otherwise it is called 
\emph{nonamenable}. 
\end{Definition} 

A variant and weaker property than the above is the following.
         
\begin{Definition}
\label{Defi:Anchored-Nonamenable}
Let $G := \left(\VV, \EE\right)$ be an infinite, locally finite, connected
graph. We define the \emph{anchored Cheeger constant} of $G$ 
with respect to the vertex $v_0$ by
\begin{equation}
\kappa'\left(G, v_0\right) := \inf \, 
\left\{ \left. \frac{\left\vert \partial W \right\vert}
                    {\left\vert W \right\vert}
\,\right\vert\, 
v_0 \in W \subseteq \VV, W \mbox{\ connected and\ } 
\left\vert W  \right\vert  < \infty \, \right\}
\label{Equ:Anchored-Cheeger-Const}
\end{equation}
where $\partial W$ is defined as above. 
The graph $G$ is said to be \emph{strongly amenable} if
$\kappa'\left(G, v_0\right) = 0$, otherwise it is called 
\emph{weakly nonamenable}. 
\end{Definition}
It is easily argued that for a connected graph $G$,
$\kappa'\left(G, v_0\right) = 0$ implies that $\kappa'\left(G, v\right) = 0$
for every vertex $v$ and so the definition of strong amenability
(or weak nonamenability) does not depend on the choice of the
vertex $v_0$. Of course, the value of the constant
$k'\left(G, v_0\right)$ may depend on the choice of $v_0$ in the
weakly nonamenable case. 
It follows by definition that 
$\kappa\left(G\right) \leq \kappa'\left(G, v_0\right)$ for any $v_0$ and so
strong amenability implies amenability. On the other hand, it is easy to show
that the two notions are not equivalent although if $G$ is transitive then
they are equivalent.

A special class of transitive graphs which are associated with 
finitely generated groups are the so-called \emph{Cayley graphs}.
\begin{Definition}
\label{Defi:Cayley-Graph}
Given a finitely generated group $\bar{G}$ and a symmetric
generating set $\SSS$ (symmetric meaning that $\SSS$=$\SSS^{-1}$),
a graph $G := \left(\VV, \EE\right)$ is called the
\emph{left-Cayley graph} of $\bar{G}$ obtained using $\SSS$ if
the vertex set of $G$ is $\bar{G}$ and the edge set is 
$\left\{ \left\{ u, v \right\} \,\Big\vert\, v = s u 
\mbox{\ for some\ } s \in \SSS \,\right\}$. 
\end{Definition}
Similarly we can also define a \emph{right-Cayley graph} of the group
$\bar{G}$ obtained using $\SSS$. Observe that 
the left- and right-Cayley graphs obtained using the same 
symmetric generating set are isomorphic, where
an isomorphism is given by the group involution 
$u \mapsto u^{-1}$, $u \in \bar{G}$. It is also easy to see that multiplication
on the right by any element in $\bar{G}$ is a graph automorphism of any 
left-Cayley graph.

If not explicitly mentioned otherwise,
by a \emph{Cayley graph} of a finitely generated group $\bar{G}$, we will 
always mean a left-Cayley graph with respect to some 
symmetric generating set.  

\begin{Definition}
A group is \emph{finitely presented} if it is described by a finite number
of generators and relations.
\end{Definition}

\begin{Definition}
A graph is \emph{one-ended} if when one removes any finite subset of the 
vertices, there remains only one infinite component. A group is
\emph{one-ended} if its Cayley graph is; it can be shown that this
is then independent of the generators used to construct the Cayley graph.
\end{Definition}

\section{Limit of the tail of the
cluster size distribution for Cayley graphs}
\label{Sec:Ordered-Group}

In this section, we prove Theorem \ref{Thm:decay}.
Throughout this section $o$ will denote the identity element of our group.

A Cayley graph is said to have \emph{polynomial growth} 
if the size of a ball is bounded by some polynomial (in its radius). 
Given a finitely generated group, its Cayley graph
having polynomial growth does not depend on the choice of the finite symmetric
generating set. It is well known (see \cite{Harpe})
that the growth of a Cayley graph of polynomial growth
is always between $\alpha \, r^k$ and ${1 \over \alpha} r^k$, 
for some $k \in \Nbold$ and $\alpha \in (0,1)$
and that if a Cayley graph is not of polynomial growth, 
then for any polynomial $p(n)$, the ball of radius $n$ around 
$o$ is larger than $p(n)$ for all but at most finitely many $n$.

Let $G$ be a Cayley graph with degree $d$. Denote by $C_x$ the open 
component of 
vertex $x$; $C$ will stand for the open component of $o$. As usual, for a 
(not necessarily induced) 
subgraph $H$ of $G$, $E(H)$ is the edge set and $V(H)$ 
is the vertex set of $H$.
Given some $p\in [0,1]$, let $\pi_n := \Pbb_p\left(|C|=n\right)$.

\begin{Lemma} 
\label{Lemma:quadratic}
If $G$ is a Cayley graph of linear or of quadratic growth, then 
\[ \lim_{n\rightarrow \infty} 
-{1\over n} \log \, \Pbb_p\left(|C|=n\right)\]
exists for every $p \in \left(0,1\right)$.  
\end{Lemma}

\proof  If $G$ has
quadratic growth then the vertices of $G$ can be partitioned into 
finite  classes, so-called {\it blocks of imprimitivity}, in 
such a way 
that the group of automorphisms restricted to the classes is $\Zbold^2$, 
see 
\cite{ST}. Now we can mimic the proof of the claim for $\Zbold^2$, see 
\cite{Gri99}: use 
the subadditive theorem and the fact that for any two connected 
finite 
subgraphs of $G$, one of them has a translate that is disjoint from the 
other, but 
at bounded distance from it.
For Cayley graphs of linear growth, one can proceed along the same
arguments, since a partition into blocks of imprimitivity, as above,
exists (see \cite{Harpe}). 
\qed 

Before starting the proof of Theorem \ref{Thm:decay},
we first prove the following lemma which 
gives an important estimate for Cayley graphs with at least cubic
growth. Using the simple structure of Cayley graphs of linear or
quadratic growth, Lemma \ref{Lemma:group} is true for every Cayley graph. 
(In the latter two cases, $(|A|+|B|)^{3/4}$ can be replaced by 1.)

\begin{Lemma}
\label{Lemma:group}
Let $G$ be a Cayley graph of at least cubic growth and $A,B\subset V(G)$. 
Then there is a $\gamma\in \Aut (G)$ such that the translate $\gamma A$ is 
disjoint from $B$ and $\dist (\gamma A, B)\leq (|A|+|B|)^{3/4}$.
\end{Lemma}

\proof
Let ${\cal A}_n$ be the set of all connected subgraphs of size $n$ in $G$ that 
contain the $o$.
Fix some group $\Gamma$ of automorphisms of $G$ such that $\Gamma$ 
is vertex-transitive on $G$ and only the identity of $\Gamma$ 
has a fixed point. 
(If $G$ is a left Cayley graph then choosing $\Gamma$ to be the group 
itself
acting with right 
multiplication works.)
For a vertex $x$ of $G$, let $\gamma_x\in\Gamma$ be the (unique) element 
of $\Gamma$ that takes $o$ to $x$. Finally, for a subgraph $H$ of $G$ 
denote by $H'$ the 1-neighborhood of $H$ (that is, the set of vertices at 
distance $\leq 1$ from $H$). Note that if $H$ is connected and $|V(H)|>1$, 
then $|V(H')|\leq d|V(H)|$, because every point of $H$ has at most $d-1$ 
neighbors outside of $H$.

Let $n,m>1$ and $A\in{\cal A}_n$, $B\in{\cal A}_m$. Suppose that for some 
$\gamma\neq\gamma' \in \Gamma$ there is a point $x$ in $A'$ such that 
$\gamma B'$ and $\gamma' B'$ both contain $x$. Then, by the choice of 
$\Gamma$, $\gamma^{-1} x\not =\gamma'^{-1}x$. Since $\gamma^{-1} x, 
\gamma'^{-1} x\in B'$, we conclude that every $x\in A'$ is contained 
in at most $|V(B')|$ translates of $B'$ by $\Gamma$. Hence there are
at most
$|V(A')|\,|V(B')|$ translates of $B'$ that intersect $A'$. Since 
$G$ has at 
least cubic growth, so there is a constant $\alpha > 0$ such that,  
the ball of radius $(n+m)^{3/4}$ 
around $o$ contains at least $\alpha (n+m)^{9/4}$ points, which is greater 
than $|V(A')|\,|V(B')|\leq d^2 nm$ for $m,n$ sufficiently large. Therefore 
there exists a vertex $x_{A,B}$ in this ball of radius $(n+m)^{3/4}$ 
such that 
$\gamma_{x_{A,B}} B'$ does not intersect $A'$. Fix such an $x_{A,B}$. 
Fix some path $P(A,B)$ of 
minimal length between $A$ and $x_{A,B}$, denote its length by $|P(A,B)|$. 
By the choice of $x_{A,B}$ we have $|P(A,B)|\leq (n+m)^{3/4}$.
Taking $\gamma := \gamma_{x_{A,B}}^{-1}$ completes the proof. \qed

\noindent
\emph{Proof of Theorem \ref{Thm:decay}:}
For graphs of linear or quadratic growth, the theorem follows from Lemma
\ref{Lemma:quadratic}.

Assume now that our group has at least cubic growth 
and so the ball of radius $r$ has volume $\geq \alpha r^3$ with some $\alpha
>0$ by the facts about Cayley graphs that we mentioned earlier.
Fix $\Gamma$ as in the proof of the previous lemma.

The generalized subadditive limit theorem 
(see Theorem II.6 in the Appendix of \cite{Gri99})
gives the result if we can show that
\begin{equation}
\label{Equ:subadditive}
 \pi_{m+n}\geq \pi_m \pi_n c^{{(m+n)}^{3/4}\logg (m+n)} 
\end{equation}
whenever $m$ and $n$ are sufficiently large, where $0<c=c(d,p)<1$
is some constant depending only on $d$ and $p$.

We will first show that
$$
2^{(m+n)^{3/4}((1+d)\logg(m+n)+c_1(d))}\pi_{m+n} 
$$
$$
 \geq  \sum_{A\in {\cal A}_n}\sum_{B\in {\cal A}_m}
           \Pbb_p\left(C=A\right)
\Pbb_p\left(C_{x_{A,B}}=\gamma_{x_{A,B}}(B)\right)
           p^{(n+m)^{3/4}}(1-p)^{2d(n+m)^{3/4}}
$$
where $c_1(d)$ is a constant depending only on $d$.
We will then show that the theorem will follow easily from here.

To prove the above inequality let $A\in{\cal A}_n$, $B\in{\cal A}_m$. 
Define $x_{A,B}$ and $\gamma_{x_{A,B}}$ as in the proof of 
Lemma \ref{Lemma:group}. Let 
$U(A,B)$ be defined as the union of three graphs:
$U(A,B):=A\cup\gamma_{x_{A,B}}B\cup P(A,B)$. Fix some arbitrary 
$\tilde{X}(A,B)\subset U(A,B)$ set of vertices not containing $o$ such that 
the subgraph
$K(A,B):=U(A,B)\setminus \tilde{X}(A,B)$ is connected and 
$|V(K(A,B))|=n+m$. Then let $X(A,B)$ be the subgraph of $U(A,B)$ 
consisting of the edges incident to some element of $\tilde{X}(A,B)$.

For fixed $A\in {\cal A}_n$ and $B\in {\cal A}_m$ we obtain
\begin{eqnarray*}
 &      & \Pbb_p\left(C=K(A,B)\right)\\
 & \geq & \Pbb_p\left(C=A\right)
          \Pbb_p\left(C_{x_{A,B}}=\gamma_{x_{A,B}} (B)\bigl| \,C=A \right)
          p^{|P(A,B)|}(1-p)^{2d|P(A,B)|} \,
\end{eqnarray*}
by first opening the edges of $P(A,B)$, closing the other edges incident 
to the 
inner vertices of $P(A,B)$ but not in $A\cup \gamma_{x_{A,B}}(B)$, 
and finally closing every edge incident to some element of $X(A,B)$, 
whenever it is necessary. 
The events $\{C=A\}$ and $\{C_{x_{A,B}}=\gamma_{x_{A,B}}(B)\}$ are 
independent because they are determined by disjoint sets of edges, since $x_{A,B}$ was chosen such that $A'$ and
$\gamma_{x_{A,B}}(B')$ are disjoint. 
Hence the previous inequality can be rewritten as
\begin{eqnarray}
 &      & \Pbb_p\left(C=K(A,B)\right) \nonumber \\ 
 & \geq & \Pbb_p\left(C=A\right)
          \Pbb_p\left(C_{x_{A,B}}=\gamma_{x_{A,B}} (B)\right)
          p^{(m+n)^{3/4}}(1-p)^{2d(m+n)^{3/4}} 
          \label{Equ:egy}
\end{eqnarray}
also using $|P(A,B)|\leq (n+m)^{3/4}$.

Now we will show that a given $K \in {\mathcal A}_{m+n}$ can be equal to
$K(A,B)$ for at most
$2^{(m+n)^{3/4}((1+d)\logg(m+n)+c_1(d))}$ pairs $(A,B)$, 
where $c_1(d)$ is a constant depending only on $d$.
First, given $m$ and $n$, $U(A,B)$
determines $(A,B)$ up to a factor 
\[2^{(\logg (m+n)+1+\logg d)(m+n)^{3/4}} m \,.\] 
This is because of the following reason.  An upper 
bound for the number of choices for the edges of $P(A,B)\setminus (A\cup 
\gamma_{x_{A,B}} B)$ from $U(A,B)$ is 
\[ |E(U(A,B))|^{(m+n)^{3/4}} \leq 
(2d(m+n))^{(m+n)^{3/4}} = 2^{(\logg (m+n)+1+\logg d)(m+n)^{3/4}} \,,
\]
using $|P(A,B)|\leq (n+m)^{3/4}$ and $|E(U(A,B))|\leq (n+m+(n+m)^{3/4})d$. 
If we delete the edges of $P(A,B)\setminus (A\cup  \gamma_{x_{A,B}} B)$ 
from 
$U(A,B)$,
we get back $A\cup \gamma_{x_{A,B}} B$. This has two components, so one of
them is $A$ and the other one is  $\gamma_{x_{A,B}} B$.
The set $\gamma_{x_{A,B}}B$ may coincide for at most $|V(B)|=m$ many 
different $B$'s (all
being $\Gamma$-translates of $\gamma_{x_{A,B}} B$ to $o$, using again the choice of $\Gamma$). We conclude 
that the number of $(A,B)$ pairs that give the same $U(A,B)$ is at most 
$2^{(m+n)^{3/4}(\logg (m+n)+1+\logg d)}m$.
Now, $X(A,B)$ is $U(A,B)\setminus E(K(A,B))$ without its isolated points 
(points of degree 0) and so for a given $K\in {\mathcal A}_{m+n}$, 
\begin{eqnarray*}
 &      & |\{(A,B): K(A,B)=K\}| \\
 & \leq & 2^{(m+n)^{3/4}(\logg (m+n)+1+\logg d)} m \, |\{U(A,B): K(A,B)=K\}| \\
 &   =  & 2^{(m+n)^{3/4}(\logg (m+n)+1+\logg d)} m \, 
          |\{X(A,B): K(A,B)=K\}|\,. 
\end{eqnarray*}

We will bound the cardinality of the set on the right side, with this 
fixed $K$.
Given $A$ and $B$, $X(A,B)$ is such a graph that the
union $K(A,B)\cup X(A,B)$ is
connected, and $|V(X(A,B))|\leq d(n+m)^{3/4}$ (since $X$ is contained in
the 1-neighborhood of $\tilde{X}$, and $|\tilde{X}|\leq (n+m)^{3/4}$). 
To find an upper bound for the
number of possible $X(A,B)$'s with these two properties (and hence where
possibly
$K(A,B)=K$), we first specify the
vertices of $K(A,B)$ that are also in $X(A,B)$
(at
most ${n+m}\choose {d(n+m)^{3/4}}$ possibilities). If $X(A,B)$ has $k$ 
components, with some arbitrary fixed ordering of 
the vertices of $G$, let $x_i$ be the first element 
of $K(A,B)\cap X(A,B)$ in the $i$'th component.
Then for each $x_i$ choose the size of the component of $X(A,B)$ that 
contains it. 
There are at most $2^{d(n+m)^{3/4}+1}$ total ways to do this because
the number of ways to express an integer $k$ as the ordered sum 
of positive integers (which would be representing the sizes 
of the different components) is at most $2^k$ and then we can sum this up
from 1 to $d(n+m)^{3/4}$ corresponding to the different possible sizes
for the vertex size of $X(A,B)$. Next, we finally choose the components 
themselves. It is known that the number of lattice animals on
$\ell$ vertices is at most $7^{2d\ell}$ (see (4.24) in \cite{Gri99})
which gives us a total bound of $7^{2d^2(n+m)^{3/4}}$
for the number of ways to choose all the components.
Note that we 
did not have to choose $x_i$, since $x_i$ is determined by $X(A,B)\cap
K(A,B)$ as soon as we know the components of the $x_j$ for all $j<i$.
Calculations similar to the above
can be found in \cite{Tim06}. We obtain an upper bound of
$$
2^{(m+n)^{3/4}(\logg (m+n)+1+\logg d)}m {{n+m}\choose
           {d(n+m)^{3/4}}}2^{d(n+m)^{3/4}+1} 7^{2d^2(n+m)^{3/4}}
$$
for the number of all possible
pairs $(A,B)$ that define the same $K=K(A,B)$ for some connected subgraph
$K$ with $n+m$ vertices, whenever $m$ and $n$ are not too small. Bounding
the binomial coefficient by $(n+m)^{d(n+m)^{3/4}}$, it easy to see that this
is at most 
$2^{(m+n)^{3/4}((1+d)\logg(m+n)+c_1(d))}$ for some constant $c_1(d)$.

Since every $K(A,B)$ is in ${\mathcal A}_{m+n}$,
the first inequality below follows from this last estimation on the
overcount. The second one is a consequence of (\ref{Equ:egy})
\begin{eqnarray*}
  &      & 2^{(m+n)^{3/4}((1+d)\logg(m+n)+c_1(d))}\pi_{m+n} \\
  & \geq & \sum_{A \in {\cal A}_n} \sum_{B \in {\cal A}_m}
           \Pbb_p \left(C=K(A,B) \right) \\
  & \geq & \sum_{A\in {\cal A}_n}\sum_{B\in {\cal A}_m}
           \Pbb_p\left(C=A\right)
\Pbb_p\left(C_{x_{A,B}}=\gamma_{x_{A,B}}(B)\right)
           p^{(n+m)^{3/4}}(1-p)^{2d(n+m)^{3/4}}\\
  & \geq & \pi_n \pi_m \beta^{(n+m)^{3/4}} \\
  & \geq & \pi_n \pi_m \beta^{(n+m)^{3/4} \logg (m+n)}\,,
\end{eqnarray*}
where $\beta := p(1-p)^{2d} \in \left(0,1\right)$,
whenever $m$ and $n$ are large enough. This yields Equation
(\ref{Equ:subadditive}) with an appropriate choice of 
$c(d,p)$ as desired and proves the theorem. \qed

\noindent
{\bf Remarks:}
The following claim seems intuitively clear, but ``continuity" arguments 
that work for $Z^d$ (or more generally, for so-called ordered groups) fail 
for arbitrary groups. If it were true, then the proof 
of Theorem \ref{Thm:decay} would become significantly simpler: the 
subadditive theorem could be applied almost right away.
\begin{Question}
\label{Q:touch}
Let $G$ be a transitive graph. Is there a constant $c$ depending on 
$G$ such that for any finite subgraphs $A$ 
and $B$ there is an automorphism $\gamma$ such that $\gamma A$ and $B$ 
are disjoint and at distance $c$ from each other?
\end{Question}

Our Lemma \ref{Lemma:group} only shows that there exists a $\gamma$ such 
that $A$ and $\gamma B$ are at distance $\leq (|A|+|B|)^{3/4}$.
As observed by Iva Koz\'akov\'a (personal communication),
one cannot have a positive answer to Question \ref{Q:touch}
with $c=1$ for all groups. 
An example showing this is the free product of a cycle 
of length 3 and a cycle of length 4.

It is worth noting that for a Cayley graph of a so-called 
\emph{ordered group}, the
proof of Theorem \ref{Thm:decay} is rather straightforward. This is primarily
because of the remarks made above. In this case the proof is really a 
generalization of the proof for $\Zbold^d$. Interesting enough one can also 
show that on the infinite regular tree with degree $3$ 
(with is not a Cayley graph of an ordered group)
such an argument does not work. Still, Theorem \ref{Thm:decay} holds 
of course for it and there is in fact
an affirmative answer to Question \ref{Q:touch} in this case.

\section{Sub-exponential decay for certain transitive amenable graphs in
the supercritical regime}
\label{Sec:Amenable-SubExp-Decay}

While Question \ref{Q:Transitive-Nonamenable-Exp-Decay} and 
Conjecture \ref{Conj:Transitive-SubExp-Decay} propose a
characterization of amenability via cluster size decay in the
supercritical regime (assuming, for completeness, the widely 
believed conjecture \cite{BenSch96}, that $p_c<1$ whenever $G$ grows 
faster than linear), a conjecture of Pete suggests that this sharp
contrast vanishes from a slightly different point of view. Instead of the
size of the cluster, consider the size of its boundary. It is known
from Kesten and Zhang \cite{Kesten-Zhang} that when $G=\Zbold ^d$, for
all $p > p_c$, there exists a $k$ such that the probability that the exterior
boundary of the $k$-closure (see Definition \ref{Defi:Closure})
of a finite supercritical cluster has size
$\geq n$ decays exponentially in $n$. (This is not true without taking the
closure, as also shown in \cite{Kesten-Zhang} for $p\in (p_c,1-p_c)$.) 
This led Pete to conjecture 
that for any transitive graph and supercritical $p$, there
exists a constant $k=k(p)$ such that $\Pbb (n<|\partial_k^+ 
C(o)|<\infty)
\leq \exp (-cn)$, where $\partial_k^+ C(o)$ denotes the exterior boundary
of the $k$-closure of the cluster of $o$. See \cite{Pete-preprint} for
applications.

Before starting on the proof of Theorem \ref{Thm:Amenable-SubExp-Decay},
we prove the following (technical) lemma which will be needed in the proof. 
\begin{Lemma}
\label{Lem:Foelner}
Let $G$ be an amenable Cayley graph. Then there is a sequence 
$\left\{W_n\right\}_{n \geq 1}$ of subsets of $\VV$
such that for every $n \geq 1$ the induced graph on $W_n$ is connected and
\begin{equation} 
\lim_{n \rightarrow \infty} \frac{\left\vert \partial W_n \right\vert}
                              {\left\vert W_n \right\vert}
= 0.
\label{Equ:Foelner-Sets}
\end{equation} 
Moreover,
$\sup_n \frac{\left\vert W_{n+1} \right\vert}{\left\vert W_n \right\vert} 
< \infty$. 
\end{Lemma}

\proof 
For groups of linear or quadratic growth, define
$W_n$ to be the ball of radius $n$ and it is immediate. (In fact,
for all groups of polynomial growth, the (nontrivial) 
facts we mentioned earlier concerning them implies that we can take
$W_n$ to be the ball of radius $n$ in these cases as well.)

We now assume that the group as at least cubic growth rate.
Since $G$ is amenable,
there exists a sequence $\left\{ W_n \right\}_{n \geq 1}$
of nonempty finite subsets of $\VV$
such that for every $n \geq 1$ the induced subgraph on $W_n$ is connected and
satisfies equation (\ref{Equ:Foelner-Sets}). 
(In the definition of amenability, the $W_n$'s are not necessarily connected, 
but it is easy to check that they may be taken to be.)
Without loss of generality, we can also assume 
$\left\vert W_n \right\vert \leq \left\vert W_{n+1} \right\vert$. 

Now, whenever $\left\vert W_{n+1} \right\vert / \left\vert W_n \right\vert > 3$
we will add a new set $E$ in the F{\o}lner sequence, after $W_n$, with the
property that $E$ is connected,
that $|E|/|W_n|\leq 3$, and such that $|\partial E|/|E|\leq
|\partial W_n|/|W_n|+2d/|W_n|^{1/4}$. 
The lemma then can be proved by repeating this procedure as long as
there are two consecutive sets in the sequence whose sizes have ratio
greater than $3$.

So all what is left, is to show the existence of such an $E$. Now, apply
Lemma \ref{Lemma:group}
with $A$ and $B$ both chosen to be $W_n$. Take the union of $A$,
$\gamma B$, and the path of length $\leq (|A|+|B|)^{3/4}$ between $A$ and
$\gamma (B)$. Let the resulting graph be $E$. Clearly $E$ satisfies the
condition about its size. It also satisfies the isoperimetric
requirement, because
$|\partial E|\leq 2|\partial W_n|+(2d |W_n|)^{3/4}$ and $|E|\geq 2|W_n|$, 
where $d$ is the degree of a vertex in $G$. 
This completes the proof. \qed

\noindent
\emph{Proof of Theorem \ref{Thm:Amenable-SubExp-Decay}: } 
Let $\left\{ W_n \right\}_{n \geq 1}$ be a sequence of subsets of $\VV$ 
satisfying the conditions of Lemma \ref{Lem:Foelner}.

For a finite set $W \subseteq V\left(G\right)$, let
$\partial_{{\tt Ext}} W$ be the set of
$v\in\partial W$ for which there exists a path from $v$ to
$\infty$ which lies (other than $v$) in $V\left(G\right) \setminus 
(W\cup \partial W)$.
It is easy to see that if the induced graph on $W$ is connected, then for any 
vertex $w \in W$ the set $\partial_{{\tt Ext}} W$ is a 
\emph{minimal cutset} between $w$ and $\infty$. 
From \cite{BabBen99, Tim06} we know that, 
since we are assuming the graph $G$ is a Cayley graph of a finitely presented group 
with one end, there exists a positive integer $t_0$, such that
any minimal cutset $\Pi$ between any vertex $v$ and $\infty$ must satisfy 
\begin{equation}
\forall A,B \mbox{ with } \Pi=A\cup B,
\mbox{dist}_G\left(A, B\right) \leq t_0.
\label{Equ:Def}
\end{equation}
Letting $U^{t} := \left\{ v \in V\left(G\right) \,\Big\vert\, 
\mbox{dist}_G\left(v, U\right) \leq t \,\right\}$
for any $U \subseteq V\left(G\right)$, and $t \in \Nbold$,
it is not hard to deduce from the above that
for any connected finite subset of vertices $W$,
we have that the induced subgraph on
$\left(\partial_{{\tt Ext}} W\right)^{t_0}$ is connected.
In particular, it follows that for each $n \geq 1$ the induced graph on
$\left(\partial_{{\tt Ext}} W_n\right)^{t_0}$ 
is connected, and further
by using (\ref{Equ:Foelner-Sets}) we get 
\begin{equation}
\lim_{n \rightarrow \infty} \frac{\left\vert 
                               \left(\partial_{{\tt Ext}} W_n \right)^{t_0}
                               \right\vert}
                              {\left\vert W_n \right\vert}
= 0 \,.
\end{equation}

Now the proof by Aizenman, Delyon and Souillard 
\cite{AizDeSou80} as given in \cite{Gri99} (see page 218), 
essentially goes through when we replace a ``$n$-ball'' of
$\Zbold^d$ by $W_n$, and the ``boundary of a $n$-ball'' by 
$\left(\partial_{{\tt Ext}} W_n\right)^{t_0}$, leading to the 
sub-exponential bound (\ref{Equ:Transitive-SubExp-Tail-Decay}). 
The point of Lemma \ref{Lem:Foelner} is that we need to obtain
the claim in the theorem for all $n$; without Lemma \ref{Lem:Foelner},
we could only make the conclusion for a sequence of $n$ going to $\infty$.
\qed

\noindent
{\bf Remarks:} Note that to carry out the above proof, we do not need that
(\ref{Equ:Def}) holds for all minimal cutsets but only for some fixed 
F{\o}lner
sequence, i.e. for a sequence of connected $W_n$'s satisfying
(\ref{Equ:Foelner-Sets}). Thus a positive answer to the
following question would imply Theorem \ref{Thm:Amenable-SubExp-Decay}
for an arbitrary amenable group. 

\begin{Definition}
\label{Defi:Closure}
The \emph{$k$-closure} of a graph $G$ is defined to 
be the graph on the vertex set of $G$ with an edge between two vertices if 
and only if their distance in $G$ is at most $k$. 
\end{Definition}

Hence
(\ref{Equ:Def}) is equivalent to saying that any minimal cutset $\Pi$ of 
$G$ is 
connected in the $t_0$-closure of $G$.

\begin{Question}
\label{Q:Foelner}
Does every amenable graph have a F{\o}lner sequence $\{W_n\}$ (that is, a
sequence satisfying Equation (\ref{Equ:Foelner-Sets})), 
such that for some $k$ the $k$-closure of
$\partial W_n$ is connected for every $n$? 
\end{Question}
In \cite{Tim06}, an example of a Cayley graph 
(coming from the so-called the lamplighter group) with one end, whose ``usual''
F{\o}lner sequence does not satisfy
the above property for its minimal cutsets, is given.
However, 
the lamplighter group is not a counterexample to Question \ref{Q:Foelner}
as shown by the following construction:

\begin{Example}
\label{Ex:Foelner}
\rm{
Recall that, informally, the lamplighter group $G$ is defined as follows.
An element of $G$ is a labeling of $\Zbold$ with labels ``on" or
``off", with only finitely many on, together with one specified element
of $\Zbold$, the
position of the lamplighter. Take the element when we move the lamplighter 
one
step to the right (corresponding to multiplication from the right by the element with 
all the lamps off and the lamplighter at 1), and the element when we 
switch 
the lamp where the
lamplighter is (corresponding to the element when the lamplighter is in 0 
and the lamp there is the only one on), as a set of generators for the 
right-Cayley graph that we 
consider now. This way we defined multiplication for any two elements. See
e.g. \cite{Tim06} for a more formal definition.

Given $x\in G$, let $\pi (x)$ be the position of the lamplighter.

To construct the desired F{\o}lner sequence $W_n$, let $B_n$ be the set of
elements $x$ with $\pi (x)\in [1,n]$ and all the lamps outside 
$[1,n]$
are off. We shall add paths to $B_n$ to get $W_n$, in the following 
way. 
For each element $x$ of the inner boundary of $B_n$ we will define a path 
$P_x$. Note that since $x$ is on the boundary of $B_n$, $\pi (x)$ is 
either $1$ or $n$. If $\pi (x)=n$, and if 
the rightmost lamp that is on is at place $n-k$, $P_x$ will be 
the following. Start from $x$, then the lamplighter moves to the 
$n+k+1$'th 
place,
switches the lamp on there, then moves back to place $n-k$, switch the 
lamp,
then move to $n+k+1$ again and switch, and then move back to $n$. The 
endpoints of $P_x$ are in the boundary of $B_n$ and the interior of $P_x$ 
is disjoint from 
$B_n$. For 
those $x$, where $\pi (x)= 1$, use the above definition but ``reflected". 
Finally,
define a path from the point where all lamps are off and the lamplighter is
in $n$ to the one where all lamps are off and he is in 1, by sending the
lamplighter to $2n$, switch, go to $-n$, switch, back to $2n$, switch, 
back to $-n$,
switch, and then go to 1. Define $W_n$ as the union of $B_n$ and all the 
$P_x$, where $x$ is some boundary point of $B_n$.

We only sketch the proof of that $W_n$ is a F{\o}lner sequence with  
boundaries that are connected in the 2-closure. We leave it for the 
interested reader to fill out the 
details.

Look at the 2-closure $G_2$ of $G$. Define a graph on the 
connected 
components
of the inner boundary $\partial B_n$ in this thickening: put an edge 
between two if some 
points
of the two are connected by a path $P_x$ defined above. One can show that 
the
paths were defined so that this graph is connected. The 
boundary of a path 
is clearly connected in the 2-closure, and (one can show that) these
path-boundaries are (basically) contained in $\partial W_n$. One concludes
that $W_n$ has a connected boundary in $G_2$.
$W_n$ is F{\o}lner, because the paths added (and hence their boundaries, 
and
the boundary of $W_n$) were constructed so as to have total length 
constant
times $2^n$, while $W_n$ has size of order $n2^n$.
}
\end{Example}

We complete this section with the 

\medskip\noindent
\emph{Proof of Proposition \ref{prop:transneeded}:}

Let $\Tbold_2^{\rho}$ be the infinite rooted binary tree, with root 
$\rho$. Consider the graph $G$ 
obtained by attaching an infinite \emph{ray} (a copy of $\Zbold_+$)
at the vertex $\rho$ of the
graph $\Tbold_2^{\rho}$. It is easy to see that $G$ is amenable even though 
$\Tbold_2^{\rho}$ is not, but of course is not transitive. 

Observe that for the i.i.d.\ bond percolation on $G$, the critical 
probability $p_c\left(G\right) = p_c\left(\Tbold_2\right) = \sfrac{1}{2}$. 
Let $C$ be the open connected component containing the vertex $\rho$.
Fix $0 < p < 1$, it is immediate that under the measure $\Pbb_p$ we have
\[
\left\vert C \right\vert \ed X + Y \,,
\]
where $X \sim \mbox{Geometric}\left(p\right)$ and
$Y$ has a distribution same as the total size of a Galton-Watson
branching process with progeny distribution 
$\mbox{Binomial}\left(2, p\right)$, and $X$ and $Y$ are independent. 

Now for every $p \in \left(0,1\right)$ we know that $X$ has an
exponential tail. Moreover from classical branching process theory
we know that
when $p \neq \sfrac{1}{2}$ we must also have an exponential tail for $Y$  
on the event $\left[ Y < \infty\right]$. This is because when 
$p < \sfrac{1}{2}$ the process is subcritical and we can use
Lemma \ref{Lem:Size-of-Sub-Critical-GW}(a) 
(given in Section \ref{Sec:Special-Non-Transitive-Graph}); and when 
$p > \sfrac{1}{2}$, the process is
supercritical, but on the event $\left[ Y < \infty \right]$ it is
distributed according to a subcritical Galton-Watson process
(see Theorem I.D.3 on page 52 of \cite{AthNey72}). 
These facts together 
prove that for every $p \neq \sfrac{1}{2}$ there is a constant
$\lambda\left(p\right) > 0$, which may depend on $p$, such that 
\begin{equation}
\Pbb_p\left( \left\vert C \right\vert = n \right) 
\leq \exp\left( - \lambda\left(p\right) n \right)
\,\,\,\, \forall \,\,\, n \geq 1 \,.
\end{equation}
\qed

\section{A special non-transitive,  nonamenable graph}
\label{Sec:Special-Non-Transitive-Graph}

In this section we will study a particular non-transitive graph, which is
also nonamenable; namely, the
$d$-dimensional integer lattice $\Zbold^d$ with rooted regular trees planted
at each vertex of it. More precisely, for each $\bx \in \Zbold^d$, 
let $\Tbold_r^{\bx}$ be an infinite
rooted regular tree with degree $(r+1)$ which is rooted at $\bx$. Thus each
vertex of $\Tbold_r^{\bx}$ has degree $(r+1)$ except for the root $\bx$, which 
has degree $r$. We consider the graph
\begin{equation}
G := \Zbold^d \bigcup \left( \bigcup_{\bx \in \Zbold^d} \Tbold_r^{\bx} \right).
\label{Equ:Lattice-with-Planted-Trees}
\end{equation}
Observe that $p_c(G)=\min\left\{p_c\left(\Zbold^d\right), 
p_c\left(\Tbold_r^{\bzero}\right)\right\}$ and recall
$p_c\left(\Tbold_r^{\bzero}\right) = \sfrac{1}{r}$.

Before we prove Theorem \ref{Thm:Non-Transitive-Two-Regime},
we note that for this particular graph $G$ the
two critical points $p_c\left(\Zbold^d\right)$ and $\sfrac{1}{r}$ play 
two different roles. $p_c\left(\Zbold^d\right)$ is the critical point for
the i.i.d.\ Bernoulli bond percolation on $G$, but when $p$ is between
$p_c\left(\Zbold^d\right)$ and $\sfrac{1}{r}$ the cluster size of the 
origin behaves like the cluster size of the origin for supercritical
bond percolation
on $\Zbold^d$. This is of course intuitively clear, because in this region 
the planted trees are all subcritical. On the other hand when 
$p > \sfrac{1}{r}$ then the tree components take over and we have the
exponential decay of the cluster size of the origin, conditioned to be
finite (see Lemma \ref{Lem:Size-of-Sub-Critical-GW} below). 

The following lemma will be needed to prove Theorem 
\ref{Thm:Non-Transitive-Two-Regime}.
This result is classical in the branching process literature, for 
a proof see \cite{Kol86, Al91}. 

\begin{Definition}
we say that a nonnegative random variable has an {\it exponential tail}
if there exists $c$ so that $P(X\ge t)\le e^{-ct}$ for $t\ge 1$.
\end{Definition}

\begin{Lemma}
\label{Lem:Size-of-Sub-Critical-GW}
Consider a subcritical or critical branching process with
progeny distribution $N$. Let $S$ be the total size of the population
starting with one individual. 
\begin{itemize}
\item[(a)] If $N$ has an exponential tail and the process is subcritical,
then $S$ has an exponential tail.

\item[(b)] If the process is critical and $N$ has a finite variance, then
then there is a constant $c$ (depending on the distribution of $N$), such that 
           \begin{equation}
           \bP\left( S = n \right) \le \frac{c}{n^{3/2}} \,
           \label{Equ:Critical-GW-Size-Asymtotic}
           \end{equation}
and if the offspring distribution is \emph{non-lattice} then
           \begin{equation}
           \bP\left( S = n \right) \sim \frac{c}{n^{3/2}} \,.
           \end{equation}
\end{itemize} 
\end{Lemma}

\noindent
\emph{Proof of Theorem \ref{Thm:Non-Transitive-Two-Regime}:} 
(a) Let $C_{\Zbold^d} := C \cap \Zbold^d$ and
$C_{\Tbold_r^{\bx}} :=  C \cap \Tbold_r^{\bx}$ for $\bx \in \Zbold^d$. 
By definition 
\begin{equation}
C = \bigcup_{\bx \in C_{\Zbold^d}} C_{\Tbold_r^{\bx}} ,
\label{Equ:Non-Transitive-Component-Decomp}
\end{equation}
where the union is a disjoint union. Thus 
\begin{equation}
\left\vert C \right\vert 
= 
\sum_{\bx \in C_{\Zbold^d}} \left\vert C_{\Tbold_r^{\bx}} \right\vert .
\label{Equ:Non-Transitive-Size-Decomp}
\end{equation}

First take
$p < p_c\left(\Zbold^d\right)$ in which case 
$\Pbb_p\left( \vert C \vert < \infty \right) = 1$. 
Using the special structure of this particular graph $G$, we observe that 
when conditioned on the random cluster $C_{\Zbold^d}$, the 
\emph{tree-components}, 
$\left\{ C_{\Tbold_r^{\bx}} \right\}_{\bx \in C_{\Zbold^d}},$ 
are independent and identically distributed, 
each being a family tree of a subcritical Galton-Watson branching
process with progeny distribution $\mbox{Binomial}\left(r, p\right)$. 
So using the representation 
(\ref{Equ:Non-Transitive-Size-Decomp}) we conclude that 
\begin{equation}
\left\vert C \right\vert \ed \sum_{j=1}^N S_j \,,
\label{Equ:Non-Transitive-Size-Decomp-dist-1}
\end{equation}
where $N \ed \left\vert C_{\Zbold^d} \right\vert$ and
$\left(S_1, S_2, \ldots \right)$ are i.i.d.\ random variables each 
distributed according to the total size of a subcritical 
Galton-Watson branching process with progeny distribution
$\mbox{Binomial}\left(r, p\right)$ and are independent of $N$. 
Now using Lemma \ref{Lem:Size-of-Sub-Critical-GW} we obtain that each 
$S_j$ has an exponential tail, and moreover from Theorem 6.75 of
Grimmett \cite{Gri99} we know that the random variable $N$ has an 
exponential tail. Then using the following (easy) lemma, who's proof is
given later, we conclude that the size of the cluster $C$ also has
an exponential tail, proving the required upper bound. 
\begin{Lemma}
\label{Lem:Exp-Tail-to-Exp-Tail}
Let $\left(S_j\right)_{j \geq 1}$ be i.i.d.\ non-negative random variables,
which are independent of $N$, which is a positive integer valued 
random variable. Let 
\begin{equation}
Z := \sum_{j=1}^N S_j \,.
\label{Equ:Defi-Random-Index-Sum}
\end{equation}
Then if $S_j$'s and $N$ have exponential tails, so does $Z$. 
\end{Lemma}

The exponential lower bound holds trivially by observing that \\
$\Pbb_p\left( n \leq \vert C \vert < \infty \right) 
\geq \Pbb_p\left( \left\vert C_{\Zbold^d} \right\vert \geq n \right)$. 

Next we take $p > \sfrac{1}{r}$, in which case
$\Pbb_p\left( \vert C \vert < \infty \right) < 1$.
Observe that 
\begin{equation}
\left[ \vert C \vert < \infty \right]
= 
\left[ \left\vert C_{\Zbold^d} \right\vert < \infty, 
\mbox{\ and\ } 
\left\vert C_{\Tbold_r^{\bx}} \right\vert < \infty 
\,\,\, \forall \,\,\, \bx \in C_{\Zbold^d} \,\right] \, .
\label{Equ:Finite-Cluster-Decomposition-Two-Regime}
\end{equation}

Once again using the special structure of this particular graph $G$, 
we observe that when conditioned on the random cluster $C_{\Zbold^d}$
and the event $\left[ \left\vert C \right\vert < \infty \right]$, the
\emph{tree-components}
$\left\{ C_{\Tbold_r^{\bx}} \right\}_{\bx \in C_{\Zbold^d}}$ 
are independent and identically distributed with
distribution being the same as a supercritical Galton-Watson tree
with progeny distribution $\mbox{Binomial}\left(r, p\right)$ conditioned
to be finite. 
%So it follows that under the conditional measure 
%$\Pbb_p \left( \cdot \,\Big\vert\, \left\vert C \right\vert  < \infty \right)$
%they are also independent of the random cluster $C_{\Zbold^d}$. 
As discussed earlier in the proof of 
Proposition \ref{prop:transneeded},
it is known that any supercritical Galton-Watson
branching process, conditioned to be finite, has the same
distribution as that 
of a subcritical Galton-Watson branching process
(see Theorem I.D.3 of \cite{AthNey72}). 
So again by Lemma \ref{Lem:Size-of-Sub-Critical-GW} the cluster sizes 
$\left\{\left\vert C_{\Tbold_r^{\bx}} \right\vert\right\}_{\bx \in C_{\Zbold^d}}$,
conditioned on $C_{\Zbold^d}$ and on $|C_{\Zbold^d}|<\infty$,
have an exponential tail.
%, under the conditional measure
%$\Pbb_p\left( \cdot \,\Big\vert\, \vert C \vert < \infty \right)$. 

Now we will show that $\left\vert C_{\Zbold^d} \right\vert$ also has an
exponential tail under the conditional measure 
$\Pbb_p\left( \cdot \,\Big\vert\, \vert C \vert < \infty \right)$. For
that we observe
\begin{eqnarray*}
\Pbb_p\left( \left\vert C_{\Zbold^d} \right\vert \geq n \,\Big\vert\,
\vert C \vert < \infty \right) 
 & = & \frac{\Pbb_p\left( \infty > \left\vert C_{\Zbold^d} \right\vert \geq n
       \,\,\, \mbox{and} \,\,\, 
       \left\vert C_{\Tbold_r^{\bx}} \right\vert < \infty \,\,\, \forall \,\,\,
       \bx \in C_{\Zbold^d} \right)}{\Pbb_p\left( \vert C \vert < \infty \right)}
       \\
 & \leq & \frac{\left(\Pbb_p\left( 
          \left\vert C_{\Tbold_r^{\bzero}} \right\vert < \infty
          \right) \right)^n}{\Pbb_p\left( \vert C \vert < \infty \right)} .
\end{eqnarray*}
Since
$\Pbb_p\left( \left\vert C_{\Tbold_r^{\bzero}} \right\vert < \infty \right) < 1$,
$\left\vert C_{\Zbold^d} \right\vert$ also has exponential tail under 
the conditional measure
$\Pbb_p\left( \cdot \,\Big\vert\, \vert C \vert < \infty \right)$.

Thus again from the representation 
(\ref{Equ:Non-Transitive-Size-Decomp}) we conclude that
under the conditional measure 
$\Pbb_p\left( \cdot \,\Big\vert\, \vert C \vert < \infty \right)$, we have
\begin{equation}
\left\vert C \right\vert \ed \sum_{i=1}^{\bar{N}} \bar{S}_i \,,
\label{Equ:Non-Transitive-Size-Decomp-dist-2}
\end{equation}
where $\left(\bar{S_1}, \bar{S_2}, \ldots \right)$ 
are i.i.d.\ random variables with
exponential tails and are independent of $\bar{N}$ which also has an
exponential tail. So finally using
Lemma \ref{Lem:Exp-Tail-to-Exp-Tail} again, 
we get the required upper bound. 

Once again the exponential lower bound can be obtained trivially by observing
\begin{eqnarray*}
\Pbb_p\left( n \leq \vert C \vert < \infty \right) 
 & \geq &  
\Pbb_p\left( n \leq \left\vert C_{\Tbold_r^{\bzero}} \right\vert < \infty, 
\mbox{\ and\ } C_{\Zbold^d} = \left\{ \bzero \right\} \right) \\
 & = & 
\left( 1- p \right)^{2d} \, 
\Pbb_p\left( n \leq \left\vert C_{\Tbold_r^{\bzero}} \right\vert < \infty \right)
\,. 
\end{eqnarray*}

\vspace{0.2in}

(b) First we obtain the lower bound for 
$p \in \left( p_c\left(\Zbold^d\right), \sfrac{1}{r} \right)$. 
We observe that 
$\Pbb_p\left( \left\vert C_{\Tbold_r^{\bx}} \right\vert < \infty \right) = 1$
for all $\bx \in \Zbold^d$. Thus from 
(\ref{Equ:Finite-Cluster-Decomposition-Two-Regime})
we conclude that under $\Pbb_p$ the events 
$\left[ \vert C \vert < \infty \right]$ and 
$\left[ \left\vert C_{\Zbold^d} \right\vert < \infty \right]$ are a.s. equal. 
So
\begin{eqnarray*}
\Pbb_p\left( n \leq \vert C \vert < \infty \right) 
 & = & \Pbb_p\left( \vert C \vert \geq n, \, 
       \left\vert C_{\Zbold^d} \right\vert < \infty \right) \\
 & \geq & \Pbb_p\left( 
          n \leq \left\vert C_{\Zbold^d} \right\vert < \infty \right) \\
 & \geq & \exp\left( - \psi_1\left(p\right) n^{(d-1)/d} \right) \, ,
\end{eqnarray*}
where $\psi_1\left(p\right) < \infty$ is a constant. The last inequality 
follows from Theorem 8.61 of Grimmett \cite{Gri99}, but as explained therein
it is much easier to derive (see page 218 of \cite{Gri99}). 

%Now to obtain an upper bound we take
%$p \in \left(p_c\left(\Zbold^d\right), \sfrac{1}{r} \right)$. 
%Once again we note that for such a value of $p$ under the measure
%$\Pbb_p$ 
Since the two events $\left[ \vert C \vert < \infty\right]$ and
$\left[ \left\vert C_{\Zbold^d} \right\vert < \infty \right]$ are
a.s. equal, we have that for every fixed $L > 0$ such that $n/L \in \Nbold$,
\begin{equation}
\Pbb_p\left( n \leq \vert C \vert < \infty \right) 
\leq \Pbb_p\left( \sfrac{n}{L} \leq 
                  \left\vert C_{\Zbold^d} \right\vert < \infty \right)
     + \bP\left( \sum_{j=1}^{n/L} S_j \geq n \right)
\label{Equ:Non-Transitive-Upper-Bound-Decomp}
\end{equation}
where $\left(S_j\right)_{j \geq 1}$ are i.i.d.\ random variables 
distributed as the total size of a subcritical Galton-Watson branching
process with $\mbox{Binomial}\left(r, p\right)$ progeny distribution. 

Now from Lemma \ref{Lem:Size-of-Sub-Critical-GW} we get that 
$\mu := \bE\left[ S_1 \right] < \infty$ and moreover the
moment generating function 
$\MMM_{S_1}\left(s\right) := \bE\left[ \exp\left(s S_1\right) \right] < \infty$
for some $s > 0$. Thus using the large deviation estimate 
Lemma 9.4 of \cite{Dur05} we will get an exponential upper bound for the 
second summand on the right hand side of 
(\ref{Equ:Non-Transitive-Upper-Bound-Decomp}), by choosing $L > \mu$. 

Moreover it follows from 
Theorem 8.65 of Grimmett \cite{Gri99} that the
first summand on the right hand side of
(\ref{Equ:Non-Transitive-Upper-Bound-Decomp}) must satisfy an
upper bound of the form
\begin{equation}
\Pbb_p\left( \sfrac{n}{L} 
             \leq \left\vert C_{\Zbold^d} \right\vert < \infty \right) 
\leq \exp\left( - \sfrac{\eta\left(p\right)}{L^{(d-1)/d}} 
                 n^{(d-1)/d} \right) \,,
\label{Equ:Grimmett-8.65}
\end{equation}
where $\eta\left(p\right) > 0$. 

This proves the required upper bound. 

\vspace{0.2in}

(c) Finally we will prove the polynomial bounds when $p = \sfrac{1}{r}$. 
First to get the upper bound, we observe that for any 
$0 < \beta < 1$ we have
\begin{equation}
\Pbb_p\left( n \leq \vert C \vert < \infty \right) 
\leq \Pbb_p\left( \lfloor n^{\beta} \rfloor \leq 
                  \left\vert C_{\Zbold^d} \right\vert < \infty \right)
     + \bP\left( \sum_{j=1}^{\lfloor n^{\beta} \rfloor} \bar{S_j} \geq n \right)
\label{Equ:Non-Transitive-Poly-Upper-Bound-Decomp}
\end{equation}
where $\left(\bar{S_j}\right)_{j \geq 1}$ are i.i.d.\ random variables
distributed as the size of a critical Galton-Watson branching
process with $\mbox{Binomial}\left(r, \sfrac{1}{r}\right)$ 
progeny distribution. 

Now consider the second summand on the right hand side of  
(\ref{Equ:Non-Transitive-Poly-Upper-Bound-Decomp}), 
\begin{eqnarray}
\bP\left( \sum_{j=1}^{\lfloor n^{\beta} \rfloor} \bar{S_j} \geq n \right) 
 & \leq & 
\bP\left( \bar{S_j} \geq \sfrac{n}{\lfloor n^{\beta} \rfloor}
\mbox{\ for some\ }  1 \leq j \leq \lfloor n^{\beta} \rfloor 
\right) \nonumber \\
 & \leq & 
n^{\beta} \, \bP\left( S_1 \geq n^{1-\beta} \right) \nonumber \\
 & \leq & n^{\beta} \, \frac{c'}{\left(n^{1-\beta}\right)^{1/2}} 
          = \frac{c'}{n^{\sfrac{1}{2} - \sfrac{3}{2}\beta}}
          \label{Equ:Non-Transitive-Poly-Upper-Bound-Second-Term} \, , 
\end{eqnarray}
where the last inequality follows from 
Lemma \ref{Lem:Size-of-Sub-Critical-GW}(b)
where $c' \equiv c'\left(\beta\right) > 0$ is a constant. 

Now once again from Theorem  8.65 of Grimmett \cite{Gri99} we get that the
first summand on the right hand side
of (\ref{Equ:Non-Transitive-Poly-Upper-Bound-Decomp}) satisfies an
upper bound of the form
\begin{equation}
\Pbb_p\left( \lfloor n^{\beta} \rfloor
             \leq \left\vert C_{\Zbold^d} \right\vert < \infty \right) 
\leq \exp\left( - \sfrac{1}{2} \eta\left(p\right) 
                 n^{\beta (d-1)/d} \right) \,,
\label{Equ:Grimmett-8.65-poly}
\end{equation}
where $\eta\left(p\right) > 0$. 
Now for fixed $\eps > 0$ we take $\beta = \sfrac{2}{3} \eps > 0$,
then the required upper bound follows using 
(\ref{Equ:Non-Transitive-Poly-Upper-Bound-Decomp}),
(\ref{Equ:Non-Transitive-Poly-Upper-Bound-Second-Term}) and 
(\ref{Equ:Grimmett-8.65-poly}).

Finally to get the lower bound, we note that as in case (b), we also have
$\left\vert C_{\Tbold_r^{\bx}} \right\vert < \infty$ a.s. with respect 
to $\Pbb_p$, for all $\bx \in \Zbold^d$ and so
%So from (\ref{Equ:Finite-Cluster-Decomposition-Two-Regime})
%we conclude that under $\Pbb_p$ the events 
$\left[ \vert C \vert < \infty \right]$ and 
$\left[ \left\vert C_{\Zbold^d} \right\vert < \infty \right]$ are a.s. equal.
Thus
\begin{eqnarray*}
\Pbb_p\left( n \leq \left\vert C \right\vert < \infty \right) 
 & =    & \Pbb_p\left( \left\vert C \right\vert \geq n, 
          \left\vert C_{\Zbold^d} \right\vert < \infty \right) \\
 & \geq & \Pbb_p\left( \left\vert C_{\Zbold^d} \right\vert < \infty
          \mbox{\ and\ }
          \left\vert C_{\Tbold_r^{\bzero}} \right\vert \geq n
          \right) \\
 &  =   & \left(1 - \theta_{\Zbold^d}\left(p\right) \right) \, 
          \Pbb_p\left( 
          \left\vert C_{\Tbold_r^{\bzero}} \right\vert \geq n \right) \\
 &  =   & \frac{c''}{n^{1/2}}
\end{eqnarray*}
where $c'' > 0$ is a constant. The last equality follows from 
Lemma \ref{Lem:Size-of-Sub-Critical-GW}(b). \qed

\noindent
{\bf Remark:} 
The above theorem does not cover the case
$p = p_c\left( \Zbold^d \right)$ and for that we would need exact tail
behavior of the cluster-size distribution for critical i.i.d.\
bond percolation on $\Zbold^d$. Unfortunately, except for $d=2$ 
(see Theorem 11.89 of \cite{Gri99}) and for large $d$ (see 
\cite{TakaGord00}), such results are largely unknown. 
\vspace{0.2in}

\noindent
We now provide a proof of Lemma \ref{Lem:Exp-Tail-to-Exp-Tail} which
is presumably well known.

\vspace{0.1in}
\noindent
\emph{Proof of Lemma \ref{Lem:Exp-Tail-to-Exp-Tail}:}
%Suppose $\bP\left( N > n \right) \leq \exp\left(- \alpha n \right)$
%for all $n > 0$, where $\alpha > 0$. So 
By assumption, there exists $\gamma > 1$ such that the
generating function 
$\phi_N\left(s\right) := \bE\left[ s^N \right] < \infty$ for all
$s < \gamma$. Similarly, there exists $c > 0$, such that
%$\bP\left( S_1 > t \right) \leq \exp\left( - c t \right)$, 
%for all $t \geq 0$. 
the moment generating function 
$\MMM_{S_1}\left(\lambda\right) 
:= \bE\left[ \exp\left( \lambda S_1 \right)\right] < \infty$ for all
$\lambda < c$. 
By the Lebesgue dominated convergence theorem,
$\MMM_{S_1}\left(\lambda\right) \rightarrow 1$
as $\lambda \downarrow 0$ and so we can find $\lambda_0 > 0$ such that 
$1 \leq \MMM_{S_1}\left(\lambda_0\right) < \gamma$. 

Now by definition (\ref{Equ:Defi-Random-Index-Sum}), the moment
generating function of $Z$ is given by 
\[
\MMM_Z\left(s\right) = \phi\left(\MMM_{S_1}\left(s\right)\right) \,. 
\]
So in particular $\MMM_Z\left(\lambda_0\right) < \infty$. 
Then by Markov inequality we get
\[
\bP\left(Z > z \right) \leq \MMM_Z\left(\lambda_0\right) \, 
                            \exp\left( - \lambda_0 z \right) \,,
\]
which completes the proof. \qed

The following theorem covers the case when 
$p_c\left(\Zbold^d\right) > \sfrac{1}{r}$, in which case 
$p_c\left(G\right) = \sfrac{1}{r}$. It is not surprising that the
intermediate regime of sub-exponential decay does not appear in this case.
\begin{Theorem}
\label{Thm:Non-transitive-One-Regime}
Suppose $p_c\left(\Zbold^d\right) > \sfrac{1}{r}$
and let $C$ be the open connected component of the origin $\bzero$ of $G$.
\begin{itemize}
\item[(a)] For $p \neq \sfrac{1}{r}$ we have
           \begin{equation}
           \exp\left(- \nu_1\left(p\right) n \right) \leq 
           \Pbb_p\left( n \leq \vert C \vert < \infty \right) 
           \leq \exp\left( - \nu_2\left(p\right) \, n \right)
           \,\,\,\, \forall n \geq 1 \, ,
           \label{Equ:Non-Transitive-Exp-Decay-One-Regime}
           \end{equation}
           where $\nu_1\left(p\right) < \infty$ and 
           $\nu_2\left(p\right) > 0$. 
\item[(b)] For $p = \sfrac{1}{r}$ the lower bound in equation 
           (\ref{Equ:Non-Transitive-Poly-Decay}) holds with the same
           constant $c_1 > 0$, and the upper bound holds for
           every $\eps > 0$ but with possibly a different constant
           $c_2' \equiv c_2'\left(\eps\right) < \infty$. 

\end{itemize}
\end{Theorem}

\proof (a) First observe that when $p < \sfrac{1}{r} = p_c\left(G\right)$ then
the same argument of the first part of the proof of 
Theorem \ref{Thm:Non-Transitive-Two-Regime}(a) applies to get the upper
and lower bounds. Moreover when 
$p \geq p_c\left(\Zbold^d\right) > \sfrac{1}{r}$, 
the tree components are in supercritical regime, 
the second part of the proof of 
Theorem \ref{Thm:Non-Transitive-Two-Regime}(a) applies for both the
upper and lower bounds. 

So all remains is to prove the exponential bounds in the intermediate case
when $\sfrac{1}{r} < p < p_c\left(\Zbold^d\right)$. The argument of the second
part of the proof of Theorem \ref{Thm:Non-Transitive-Two-Regime}(a)
again goes through here, although in this case the only real
difference is 
$\Pbb_p\left( \left\vert C_{\Zbold^d} \right\vert < \infty \right) = 1$, 
which does not affect the proof. 

\vspace{0.2in}

(b) For this we follow exactly the same steps of the proof of part (c)
of Theorem \ref{Thm:Non-Transitive-Two-Regime} for both the upper and
lower bounds. The only difference here is in equation 
(\ref{Equ:Grimmett-8.65-poly}), which in
this case should be
\begin{equation}
\Pbb_p\left( \lfloor n^{\beta} \rfloor
             \leq \left\vert C_{\Zbold^d} \right\vert < \infty \right) 
\leq \exp\left( - \sfrac{1}{2} \eta\left(p\right) 
                 n^{\beta} \right) \,,
\label{Equ:Grimmett-8.65-poly-2}
\end{equation}
where $\eta\left(p\right) > 0$. So we may need a different constant than
$c_2$ in the upper bound in
equation (\ref{Equ:Non-Transitive-Poly-Decay}). \qed

\noindent
{\bf Remark:} 
Once again, the case 
$p = p_c\left(\Zbold^d\right) = \sfrac{1}{r}$ is left open, because
of similar reason as mentioned in the remark after the proof of
Theorem \ref{Thm:Non-Transitive-Two-Regime}.  
(Of course, one would be surprised if there were
any $d$ and $r$ (other than $d=r=2$) where the above held.)

\section*{Acknowledgments}
We are grateful to David J. Aldous, Jean-Francois Le Gall, Yuval Peres 
and G\'{a}bor Pete for helpful discussions. 
Bandyopadhyay wishes to thank the Department of Mathematics, 
Chalmers University of Technology, G\"{o}teborg, Sweden for its support. 
Research partially supported by the Swedish Natural Science 
Research Council (Bandyopadhyay and Steif), the G\"{o}ran Gustafsson 
Foundation for Research in Natural Sciences and Medicine (Steif) and the
Hungarian National Foundation for Scientific Research Grant TO34475 (Tim\'ar).

\bibliographystyle{plain}

\bibliography{Cluster-Size.bib}

\def\cprime{$'$}
\begin{thebibliography}{10}

\bibitem{AB87}
Michael Aizenman and David~J. Barsky.
\newblock Sharpness of the phase transition in percolation models.
\newblock {\em Comm. Math. Phys.}, 108(3):489--526, 1987.

\bibitem{AizDeSou80}
Michael Aizenman, Fran{\c{c}}ois Delyon, and Bernard Souillard.
\newblock Lower bounds on the cluster size distribution.
\newblock {\em J. Statist. Phys.}, 23(3):267--280, 1980.

\bibitem{AN84}
Michael Aizenman and Charles~M. Newman.
\newblock Tree graph inequalities and critical behavior in percolation models.
\newblock {\em J. Statist. Phys.}, 36(1-2):107--143, 1984.

\bibitem{Al91}
David Aldous.
\newblock Asymptotic fringe distributions for general families of random trees.
\newblock {\em Ann. Appl. Probab.}, 1(2):228--266, 1991.

\bibitem{AthNey72}
Krishna~B. Athreya and Peter~E. Ney.
\newblock {\em Branching processes}.
\newblock Springer-Verlag, New York, 1972.
\newblock Die Grundlehren der mathematischen Wissenschaften, Band 196.

\bibitem{BabBen99}
Eric Babson and Itai Benjamini.
\newblock Cut sets and normed cohomology with applications to percolation.
\newblock {\em Proc. Amer. Math. Soc.}, 127(2):589--597, 1999.

\bibitem{BLSP99a}
I.~Benjamini, R.~Lyons, Y.~Peres, and O.~Schramm.
\newblock Group-invariant percolation on graphs.
\newblock {\em Geom. Funct. Anal.}, 9(1):29--66, 1999.

\bibitem{BLSP99b}
Itai Benjamini, Russell Lyons, Yuval Peres, and Oded Schramm.
\newblock Critical percolation on any nonamenable group has no infinite
  clusters.
\newblock {\em Ann. Probab.}, 27(3):1347--1356, 1999.

\bibitem{BenSch96}
Itai Benjamini and Oded Schramm.
\newblock Percolation beyond {$\bold Z\sp d$}, many questions and a few
  answers.
\newblock {\em Electron. Comm. Probab.}, 1:no.\ 8, 71--82 (electronic), 1996.

\bibitem{ChenPeres04}
Dayue Chen and Yuval Peres.
\newblock Anchored expansion, percolation and speed.
\newblock {\em Ann. Probab.}, 32(4):2978--2995, 2004.
\newblock With an appendix by G\'abor Pete.

\bibitem{Harpe}
Pierre de~la Harpe.
\newblock {\em Topics in geometric group theory}.
\newblock Chicago Lectures in Mathematics. University of Chicago Press,
  Chicago, IL, 2000.

\bibitem{Dur05}
Richard Durrett.
\newblock {\em Probability: theory and examples}.
\newblock Duxbury Press, Belmont, CA, third edition, 2005.

\bibitem{Gri99}
Geoffrey Grimmett.
\newblock {\em Percolation}, volume 321 of {\em Grundlehren der Mathematischen
  Wissenschaften}.
\newblock Springer-Verlag, Berlin, second edition, 1999.

\bibitem{Hag00}
Olle H{\"a}ggstr{\"o}m.
\newblock Markov random fields and percolation on general graphs.
\newblock {\em Adv. in Appl. Probab.}, 32(1):39--66, 2000.

\bibitem{HagPeres99}
Olle H{\"a}ggstr{\"o}m and Yuval Peres.
\newblock Monotonicity of uniqueness for percolation on {C}ayley graphs: all
  infinite clusters are born simultaneously.
\newblock {\em Probab. Theory Related Fields}, 113(2):273--285, 1999.

\bibitem{HagPeresSch99}
Olle H{\"a}ggstr{\"o}m, Yuval Peres, and Roberto~H. Schonmann.
\newblock Percolation on transitive graphs as a coalescent process: relentless
  merging followed by simultaneous uniqueness.
\newblock In {\em Perplexing problems in probability}, volume~44 of {\em Progr.
  Probab.}, pages 69--90. Birkh\"auser Boston, Boston, MA, 1999.

\bibitem{HagSchStei00}
Olle H{\"a}ggstr{\"o}m, Roberto~H. Schonmann, and Jeffrey~E. Steif.
\newblock The {I}sing model on diluted graphs and strong amenability.
\newblock {\em Ann. Probab.}, 28(3):1111--1137, 2000.

\bibitem{TakaGord00}
Takashi Hara and Gordon Slade.
\newblock The scaling limit of the incipient infinite cluster in
  high-dimensional percolation. {I}. {C}ritical exponents.
\newblock {\em J. Statist. Phys.}, 99(5-6):1075--1168, 2000.

\bibitem{Jon99}
Johan Jonasson.
\newblock The random cluster model on a general graph and a phase transition
  characterization of nonamenability.
\newblock {\em Stochastic Process. Appl.}, 79(2):335--354, 1999.

\bibitem{JonSti99}
Johan Jonasson and Jeffrey~E. Steif.
\newblock Amenability and phase transition in the {I}sing model.
\newblock {\em J. Theoret. Probab.}, 12(2):549--559, 1999.

\bibitem{Kesten-Zhang}
H.~Kesten and Y.~Zhang.
\newblock The probability of a large finite cluster in supercritical bernoulli
  percolation.
\newblock {\em Ann. Probab.}, 18:537--555, 1990.

\bibitem{Kes59}
Harry Kesten.
\newblock Symmetric random walks on groups.
\newblock {\em Trans. Amer. Math. Soc.}, 92:336--354, 1959.

\bibitem{Kol86}
Valentin~F. Kolchin.
\newblock {\em Random mappings}.
\newblock Translation Series in Mathematics and Engineering. Optimization
  Software Inc. Publications Division, New York, 1986.
\newblock Translated from the Russian, With a foreword by S. R. S. Varadhan.

\bibitem{M86}
M.~V. Men{\cprime}shikov.
\newblock Coincidence of critical points in percolation problems.
\newblock {\em Dokl. Akad. Nauk SSSR}, 288(6):1308--1311, 1986.

\bibitem{Pete-preprint}
G\'abor Pete.
\newblock Anchored isoperimetry, random walks and percolation: a survey with
  many questions.
\newblock preprint.

\bibitem{Sch01}
Roberto~H. Schonmann.
\newblock Multiplicity of phase transitions and mean-field criticality on
  highly non-amenable graphs.
\newblock {\em Comm. Math. Phys.}, 219(2):271--322, 2001.

\bibitem{ST}
Norbert Seifter and Vladimir Trofimov.
\newblock Automorphism groups of graphs with quadratic growth.

\bibitem{Tim06}
{\'A}d{\'a}m Tim{\'a}r.
\newblock Cutsets in infinite graphs.
\newblock {\em Combin. Probab. Comput.}, 16:1--8, 2006.

\end{thebibliography}

\end{document}